\documentclass{amsart}

\usepackage{amsmath}
\usepackage{amssymb}
\usepackage{hyperref}

\allowdisplaybreaks

\newtheorem{theorem}{Theorem}[section]
\newtheorem{lemma}[theorem]{Lemma}
\newtheorem{proposition}[theorem]{Proposition}

\newtheorem{corollary}[theorem]{Corollary}

\theoremstyle{definition}

\newcommand{\Z}{\mathbb{Z}}
\newcommand{\R}{\mathbb{R}}
\newcommand{\C}{\mathbb{C}}

\newcommand{\nN}{\mathbb{N}}

\DeclareMathOperator{\sech}{sech}

\let\Im=\undefined\DeclareMathOperator*{\Im}{Im}

\newcommand{\qtq}[1]{\quad\text{#1}\quad}

\newcommand{\jb}[1]{\langle #1 \rangle}

\renewcommand{\l}{\ell}

\numberwithin{equation}{section}
\numberwithin{theorem}{section}

\DeclareMathOperator{\diag}{diag}
\newcommand{\dg}[1]{\diag({#1})}
\newcommand{\ceil}[1]{\lceil {#1} \rceil}

\begin{document}

\title[Growth of Fourier--Lebesgue norms for mKdV]{Growth of Fourier--Lebesgue norms for mKdV}

\author{Saikatul Haque}
\address{
Department of Mathematics\\
University of California\\Los Angeles\\CA 90095\\USA\\
\& Harish-Chandra Research Instite, Allahabad 211019, India}
\email{saikatul@math.ucla.edu}

\author{Rowan Killip}
\address{
Department of Mathematics\\
University of California\\Los Angeles\\CA 90095\\USA\\\& CEREMADE, CNRS\\Universit\'e Paris Dauphine--PSL\\ Place du Mar\'echal de Lattre de Tassigny\\ 75016 Paris\\ France}
\email{killip@math.ucla.edu}

\author{Monica Vi\c{s}an}
\address{
Department of Mathematics\\
University of California\\Los Angeles\\CA 90095\\USA}
\email{visan@math.ucla.edu}

\author{Yunfeng Zhang}
\address{Department of Mathematical Sciences\\
University of Cincinnati\\
Cincinnati\\OH 45221\\USA}
\email{zhang8y7@ucmail.uc.edu}

\begin{abstract}
We demonstrate inflation of Fourier--Lebesgue norms for solutions to the focusing modified Korteweg--de Vries equation posed on the real line.  For $p\neq 2$ and all $s\in \R$, we construct a sequence of solutions $u_n$ whose initial data $u_n(0)$ converges to zero in the Fourier--Lebesgue spaces $\mathcal{F}L^p_s(\R)$, but whose evolutions at later times $t_n$ diverge to infinity.
\end{abstract}


\maketitle


\section{Introduction}

We consider the focusing modified Korteweg--de Vries equation
\begin{equation}\tag{mKdV}\label{mkdv}
    \partial_t u+\partial_x^3u=-6u^2\partial_xu 
\end{equation} 
posed on the real line. Solutions are real-valued functions $u:\R_t\times\R_x\to\R$. 

This model is completely integrable and its plethora of conserved quantities have played an important role in the development of its well-posedness theory, particularly globally in time.  Recent developments (discussed below) have led to significant advances in understanding how these conservation laws control key norms of interest beyond the usual Sobolev spaces $H^s(\R)$ with integer $s\geq 0$.  Our goal in this paper is to demonstrate the opposite, namely, that the conservation laws do not provide control over \emph{any} Fourier--Lebesgue norms, excepting the cases where these coincide with $H^s(\R)$ spaces.

The simplest Fourier--Lebesgue norms are just the $L^p$-norms (of Lebesgue) applied to the Fourier transform.  Our convention here is
\begin{equation}\label{FT}
\widehat f(\xi):= \tfrac{1}{\sqrt{2\pi}}\int_\R e^{-ix\xi} f(x)\, dx. 
\end{equation}
More generally, given $1\leq p<\infty$ and $s\in \R$, we define
\begin{equation}\label{FLp}
\| f \|_{\mathcal{F}L^p_s(\R)} := \|  \langle\xi\rangle^s \widehat{f}(\xi) \|_{L^p(\R)}
\end{equation}
and then write $\mathcal{F}L^p_s(\R)$ for the completion of Schwartz space under this norm.  When $p=2$, this coincides with $H^s(\R)$.  In all other cases, we have the following:

\begin{theorem}\label{T:main}
For each $1\leq p<\infty$ with $p\neq 2$ and $s\in \R$, there is a sequence of Schwartz solutions $u_n(t)$ to \eqref{mkdv} that satisfy
\begin{equation}\label{E:main}
\lim_{n\to\infty} \|u_n(0)\|_{\mathcal{F}L^p_s(\R)}=0 \qtq{and}  \lim_{n\to\infty}\sup_{t\in\R} \|u_n(t)\|_{\mathcal{F}L^p_s(\R)}=\infty.
\end{equation}
\end{theorem}

Clearly, these real-valued examples apply equally well to the more general complex modified Korteweg--de Vries equation of Hirota \cite{MR338587},
\begin{align}\tag{mKdV$_\C$}\label{HmKdV}
   \partial_t u+\partial_x^3u=\pm6|u|^2\partial_x u ,
\end{align}
yielding norm inflation for this equation as well.

It has long been known that \eqref{E:main} cannot hold when $p=2$ and $s\geq0$ is an integer.   When $s=0$, for example, we have $\mathcal{F}L^2_0(\R) = L^2(\R)$ and this norm is exactly conserved by the flow.  This forms the base step of an inductive argument introduced by Lax \cite{MR369963}, which combines the traditional polynomial conservation laws and the Gagliardo--Nirenberg inequality to show that 
\begin{equation}\label{Hs int}
  \sup_t \| u(t) \|_{H^s} \lesssim \| u(0) \|_{H^s} \bigl( 1 +  \| u(0) \|_{H^s}^2\bigr)^{s} \qquad\text{for any integer $s\geq 0$,}
\end{equation}
not only for solutions to \eqref{mkdv}, but also for solutions to \eqref{HmKdV}.

Combining the time-translation symmetry with \eqref{Hs int} then shows that the $H^s$-norm of the solution at any one time controls this norm at all other times, both from above and from below.

Evidently, \eqref{E:main} and \eqref{Hs int} represent antithetical behaviors.  In fact, the motivation for the work described herein was to determine whether some bound of the type \eqref{Hs int} could also hold for Fourier--Lebesgue norms.  This is a question we were asked directly by Zihua Guo (private communication, March 2024).  A positive answer would have been useful for turning known local well-posedness results into global results.  Moreover, as we will discuss below, there are a number of positive results that led credence to this idea.  In this way, the failure of such estimates as demonstrated by Theorem~\ref{T:main} came as quite a surprise.

The bounded/growth dichotomy for all $H^s(\R)$ norms has recently been fully resolved.   In \cite{MR3820439,MR3874652}, it was shown that
\begin{equation}\label{Hs nint}
  \sup_t \| u(t) \|_{H^s} \lesssim \| u(0) \|_{H^s} \bigl( 1 +  \| u(0) \|_{H^s}^2\bigr)^{c(s)} 
\end{equation}
holds for all solutions to \eqref{HmKdV} and every real number $s>-\frac12$.  On the other hand, such estimates cannot hold when $s\leq -\frac12$; see \cite[Appendix A]{MR4726498}.  Concretely, the $H^s$-analogue of \eqref{E:main} was shown to hold for \eqref{mkdv} when $s<-\frac12$ and for \eqref{HmKdV} when $s\leq -\frac12$.  Finally, a slightly different statement was proved for \eqref{mkdv} with $s=-\frac12$, namely, for any sequence of times $t_n\to 0$, there is a sequence of Schwartz solutions $u_n$ to \eqref{mkdv} such that
$$
\|u_n(0)\|_{\dot H^{-\frac12}}\lesssim 1 \qtq{and} \lim_{n\to \infty}\|u_n(t_n)\|_{H^{-\frac12}}= \infty.
$$
This still suffices to preclude bounds of the form \eqref{Hs nint}.

The bounds \eqref{Hs nint} apply equally to the nonlinear Schr\"odinger equation
\begin{align}\label{NLS}\tag{NLS}
i\partial_t q = - \partial_x^2 q -  2|q|^2q
\end{align}
and the restriction $s>-\frac12$ is likewise sharp; see \cite{MR4726498,MR3820439,MR3874652}.

Aside from the Sobolev spaces $H^s(\R)$, there are two main classes of spaces that have been influential in the development of the well-posedness theory for \eqref{mkdv} and related equations.  These are the Fourier--Lebesgue spaces defined in \eqref{FLp} and their close cousins, the modulation spaces $M^{s,2}_p(\R)$, whose norms are given by
\begin{equation}\label{E:mod}
\|f\|_{M^{s,2}_p}:= \bigl\| \langle k\rangle^s \| \widehat f{\mkern 2mu}\|_{L^2([k, k+1])}\bigr\|_{\ell^p_k(\Z)}.
\end{equation}
When $p=2$, the modulation norm is equivalent to that on $H^s(\R)$.  If we were to replace the inner-most $L^2$-norm (over unit intervals) by the $L^p$-norm, we would obtain a norm equivalent to that of the Fourier--Lebegue space $\mathcal F L^p_s$.  Correspondingly, there is no distinction between such modulation spaces and the Fourier--Lebesgue spaces on the torus (where Fourier space is discrete).

Building on the earlier works \cite{MR3820439,MR4534495,MR4081534,MR4235636}, we have recently shown in \cite{HKVZ:1} that
\begin{align*}
\sup_{t\in\R}\, \|u(t)\|_{M^{s,2}_p}\lesssim \bigl(1+\|u(0)\|_{M^{s,2}_p}\bigr)^{c(s,p)}\|u(0)\|_{M^{s,2}_p}
\end{align*}
whenever $1\leq p<\infty$ and $0\leq s<\frac32-\frac1p$.  In fact, \cite{HKVZ:1} also proves global well-posedness of \eqref{HmKdV} in these spaces by combining   an improvement of these bounds (that propagates equicontinuity) with the sharp $H^s(\R)$ well-posedness results of \cite{MR4726498}.

We are not aware of any global well-posedness results for \eqref{mkdv} in Fourier--Lebesgue spaces with $p\neq 2$.
Nevertheless, there has been important work on the question of local well-posedness.

The study of dispersive PDE in Fourier--Lebesgue spaces began in earnest with the paper \cite{MR1820017}, which considered nonlinear Schr\"odinger equations.  The presence of a derivative in the nonlinearity of \eqref{mkdv} makes this latter model much more difficult to study.  For a very long time, $H^s(\R)$ well-posedness for \eqref{mkdv} was confined to $s\geq\frac14$, which is $\frac34$ derivatives away from scaling.   Inspired by this, \cite{MR2096258,MR2529909} tackled the well-posedness problem in Fourier--Lebesgue spaces, proving local well-posedness in $\mathcal FL^p_s$ for $2\leq p<\infty$ and $s\geq \tfrac1{2p}$.  Excitingly, this comes arbitrarily close to scaling-criticality as $p\to\infty$. Alas, Theorem~\ref{T:main} dashes the hope of globalizing this result using conservation laws.  It also crushes the hope of adapting the method of \cite{HKVZ:1} to prove $\mathcal FL^p_s$ well-posedness.

Theorem~\ref{T:main} will be proved by a careful analysis of certain exact multisoliton solutions $u_{N,\lambda}(t,x)$.  Concretely, we rescale the solutions $u_N$ with initial data
\begin{equation}\label{sech}
u_{N}(0,x) = (-1)^N N \sech( x), \qquad N\in\nN, 
\end{equation}
according to the scaling symmetry of \eqref{mkdv}:
\begin{equation}\label{resech}
u_{N,\lambda}(t,x)  := \lambda^{-1} u_N\bigl(\lambda^{-3}t,\lambda ^{-1}x\bigr), \qquad \lambda>0 .
\end{equation}

The sign choice here is of no great significance because of the $u\mapsto -u$ symmetry of \eqref{mkdv}; however, the factor $(-1)^N$ turns out to be convenient for simplifying intermediate computations.

The initial data \eqref{sech} is exactly $\pm N$ times the traditional soliton associated with \eqref{mkdv}.  Satsuma--Yajima \cite{MR463733} discovered that this is a multisoliton solution precisely when $N$ is an integer and determined the relevant spectral parameters.  

As an exact multisoliton, the evolution of $u_N$ has been understood (at least formally) for a long time: it fissions into $N$ individual solitons with characteristic widths $1,\ \tfrac13,\ \tfrac15,\ \ldots,\ \tfrac1{2N-1}$, which proceed to separate along linear trajectories as time progresses.  The asymptotic speeds of these solitons are $1,3^2,5^2,\ldots,(2N-1)^2$.  The fact that these are all squares will be technically very convenient later.   

The fact that the initial data \eqref{sech} of unit width evolves to produce a solution with much narrower solitons is eye-catching.  Taking the Fourier transform, this becomes a spectral broadening phenomenon.   This spectral broadening also occurs for \eqref{NLS} and is of great utility in optics, where it is known as super\-continuum generation; see \cite{Agrawal,MR2319053}.

The solutions $u_n$ employed in Theorem~\ref{T:main} correspond to choosing sequences of integers $N$ and scaling parameters $\lambda$ that both go to infinity in a coordinated way (depending on $p$).  In all cases, however, $\lambda \gg N$ and correspondingly, the solution $u_{N,\lambda}$ remains at low frequencies for all time.  In fact, we will provide good control on the high-frequency tails and correspondingly, the parameter $s$ has little effect.

Nevertheless, the solutions $u_{N,\lambda}$ do experience substantial broadening: the characteristic frequency scale grows from $\lambda^{-1}$ at the initial time to $N\lambda^{-1}$ at very large times.  This broadening of the Fourier spectrum is ultimately responsible for Theorem~\ref{T:main}, although it will take some effort to tease this out.  First of all, $u_{N,\lambda}$ fissions into a great many wavepackets of varying widths and one must understand the full combination.  It is helpful that these wavepackets are widely separated in physical space and so their Fourier transforms admit very different modulations.  This provides a mild form of decoupling.

In fact, we are able to exploit a stronger form of decoupling by averaging in time.  This takes advantage of the fact that the speeds are all squares and Bourgain~\cite{MR1029904} provides good $L^p$-bounds for the resulting sums; see Lemma~\ref{L:Bourgain}.   When $2<p<\infty$, we use this technique to show that the average $\mathcal FL^p_s$ norm of $u_{N,\lambda}$ is small over time intervals $[T,2T]$ with $T$ sufficiently large.  Evidently, this guarantees that there is at least one time when the norm is small.  The proof of Theorem~\ref{T:main} is then completed by employing the time-reversal symmetry. 

In the case $1\leq p<2$, we employ a duality argument and the smallness just described to show that $u_{N,\lambda}$ has large $\mathcal FL^p_s$ norm at some large time.  In this case we do not need to use the time-reversal symmetry.


Given that the proof of Theorem~\ref{T:main} hinges crucially on the exact form of the long-time asymptotics of $u_{N}$, we felt that it was important to provide the complete details of this derivation with a bare minimum of prerequisites.  For this reason, we perform a direct verification of the explicit formula for multisolitons in Section~\ref{S:2}; see Proposition~\ref{P:psi}.  In our experience,  textbook presentations of this material usually focus on NLS, rather than \eqref{HmKdV}, and employ all the machinery of inverse scattering (en route to higher goals).

In Section~\ref{S:3} we provide an elementary presentation of the long-time asymptotics of our particular solutions $u_N(t,x)$.  Due to the structure of the explicit formula \eqref{E:P:psi}, this is a matter of analyzing the asymptotics of the inverse of a concrete matrix.  There is a competing description in terms of determinants (which is equivalent via Cramer's rule); however, we do not find this to be conducive to the analysis.   Instead, we employ a well-chosen conjugation and a variation of the Schur-complement idea that partitions the matrix into nine blocks; see \eqref{A-1}.  

Finally, in Section~\ref{S:4} we analyze the Fourier--Lebesgue norms of $u_{N,\lambda}$ for large times and complete the proof of Theorem~\ref{T:main}.

\subsection*{Acknowledgements} S.H. was supported by a Fulbright-Nehru Fellowship and DST-INSPIRE grant 04/2022/001457. R.~K. was supported by NSF grants DMS-2154022 and DMS-2452346. M.~V. was supported by NSF grant DMS-2348018 and the Simons Fellowship SFI-MPS-SFM-00006244. We are grateful to Zihua Guo and Kenji Nakanishi for bringing this problem to our attention.

\section{Multisoliton solutions}\label{S:2}

The main purpose of this section is to provide a direct verification of a well-known formula for multisoliton solutions to \eqref{mkdv} using only linear algebra.   As we may do so efficiently, we present a more general formula capturing the first four flows in the \eqref{mkdv} hierarchy.

We start by introducing the key objects.  We fix $N\geq 1$, which encodes the complexity of the multisoliton. For each $1\leq j\leq N$, let $\lambda_j, a_j\in \C$ with $\Im \lambda_j>0$.  We define $\C^N$-valued functions $e$ and $\gamma$ of the spatial variable $x\in \R$ and four flow parameters $\theta, y,s,t\in\R$ as follows:
\begin{align}\label{e, gamma}
e_j = 1  \qtq{and}  \gamma_j = a_j \exp\bigl(i\theta+ i\lambda_j (x-y)  - i \lambda_j^2 s + i\lambda_j^3 t\bigr),
\end{align}
as well as $N\times N$ matrices with entries
\begin{align}\label{M}
M_{jk} = i \frac{1+\gamma_j\overline{\gamma}_k}{\lambda_j - \overline{\lambda}_k}  \qtq{and}  \Lambda_{jk} = \lambda_j \delta_{jk}.
\end{align}
Observe that
\begin{align}\label{66}
M' = - \gamma\gamma^\dagger   \qtq{and} \Lambda M - M \Lambda^\dagger  = ie e^\dagger  +i \gamma\gamma^\dagger  = ie e^\dagger  - i M',
\end{align}
where primes indicate derivatives in $x$ and ${}^\dagger $ represents conjugate transpose.

\begin{proposition}[Multisolitons]\label{P:psi}
Let $e, \gamma, M$ be as defined in \eqref{e, gamma} and \eqref{M}.  Then the function
\begin{align}\label{E:P:psi}
u:\R^4\times\R\to \C \qtq{defined by} u(\theta, y,s,t,x) = e^\dagger  M^{-1} \gamma
\end{align}
satisfies the following evolution equations:
\begin{alignat*}{2}
u_\theta &= i u,  \qquad\qquad & u_y &= -u' , \\
  iu_s&=- u'' - 2|u|^2u, \qquad\qquad  & u_t &= - u''' - 6|u|^2u',
\end{alignat*}
where the subscripts indicate partial derivatives.  These equations encode phase rotation, translation, \eqref{NLS}, and \eqref{mkdv}, respectively.
\end{proposition}

The proof of Proposition~\ref{P:psi} and the asymptotics of such solutions discussed in the next section rely crucially on several identities for Cauchy matrices.  These can be found, for example, in \cite{MR105798}.

\begin{lemma}[Cauchy Matrices]\label{L:Cauchy}
Fix $N\geq 1$ and suppose $a, b \in \C^N$ satisfy $a_j+b_k\neq 0$ for all $1\leq j, k\leq N$. Then the  $N\times N$ matrix with entries
\begin{align*}
 C_{jk} = \frac{1}{a_j +b_k}
\end{align*}
has determinant
\begin{align*}
\det C = \biggl[ \prod_{j<k}  (a_k-&a_j)(b_k - b_j) \biggr] \bigg/ \biggl[ \prod_{j,k}  (a_j + b_k) \biggr].
\end{align*}
Moreover, if $a_j\neq a_k$ and $b_j\neq b_k$ for all $j\neq k$, then the matrix $C$ is invertible and 
\begin{align*}
 ( C^{-1} )_{jj} = (a_j +& b_j)  \prod_{k\neq j}\frac{(a_k+b_j)(b_k+a_j)}{(a_k-a_j)(b_k-b_j)}, \\
 ( C^{-1} e )_j = (a_j+b_j)  \prod_{k\neq j} \frac{b_j+a_k}{b_j-b_k} &\qtq{and} (e^\dagger  C^{-1})_j= (a_j+b_j)  \prod_{k\neq j}\frac{a_j+b_k}{a_j-a_k} .
\end{align*}
\end{lemma}

The following is a direct corollary of Lemma~\ref{L:Cauchy}:

\begin{corollary}\label{C:Cauchy} Fix $N\geq 1$ and $\lambda\in \C^N$ with $\Im \lambda_j>0$ for each $1\leq j\leq N$. The Hermitian matrix with entries
\begin{align}
A_{jk} = \frac{i}{\lambda_j - \overline{\lambda}_k} \qquad \text{has} \qquad 
 \det A =  \prod_{j<k} \frac{|\lambda_j-\lambda_k|^2}{|\lambda_j-\overline{\lambda}_k|^2} \cdot \prod_{j=1}^N \frac{1}{2\Im \lambda_j}.
\end{align}
Moreover, $A$ is a positive definite matrix and so is the matrix $M$ defined in \eqref{M}.
\end{corollary}
 
\begin{proof}
The formula for the determinant of $A$ follows from Lemma~\ref{L:Cauchy}.  This formula applies also to all principal minors of $A$.  Thus, the assertion that $A$ is positive definite follows from Sylvester's criterion.

Finally, as $A$ is positive definite, for any $v\in\C^N\setminus\{0\}$ we have
$$
v^\dagger  M v = \sum_{1\leq j, k\leq N}  \overline v_j A_{jk} v_k + \overline v_j \gamma_j A_{jk} v_k \overline \gamma_k  >0,
$$
which proves that $M$ is positive definite, and so invertible.
\end{proof}

\begin{proof}[Proof of Proposition~\ref{P:psi}]
It is evident from the definitions that $M$ is independent of $\theta$,  while $\gamma$ varies according to $\partial_\theta \gamma = i\gamma$.  Thus $u_\theta = i u$.  It is also evident that $y$ acts solely as a spatial translation in $\gamma$, and so also in $M$.  Thus $u_y = -u'$.   We now move on to the more serious assertions concerning \eqref{NLS} and \eqref{mkdv}.

Straightforward computations yield
\begin{gather}
\gamma' = i \Lambda \gamma, \quad\qquad\qquad \gamma'' = - \Lambda^2 \gamma, \qquad\qquad\quad  \gamma''' = - i \Lambda^3 \gamma \label{gamma'},\\
M ' \!= - \gamma\gamma^\dagger , \quad M''\!= i\Lambda M'\! - i M' \Lambda^\dagger , \quad M''' \!= - \Lambda^2 M' \!+ 2 \Lambda  M' \Lambda^\dagger  \!- \!M' (\Lambda^\dagger )^2 \label{M'},\\
i \gamma_s = \Lambda^2 \gamma = - \gamma'' \qquad\text{and}\qquad \gamma_t = i\Lambda^3 \gamma = - \gamma'''\label{gamma_s gamma_t}, \\
iM_s = - i\Lambda M' - iM' \Lambda^\dagger  \qquad\text{and}\qquad M_t = \Lambda^2 M' + \Lambda  M' \Lambda^\dagger  + M' (\Lambda^\dagger )^2 \label{M dot},
\end{gather}
as well as
\begin{align}
u'  &= e^\dagger  M^{-1} \gamma' - e^\dagger  M^{-1} M'  M^{-1} \gamma,		\label{psi'}\\
u'' &= e^\dagger  M^{-1} \gamma'' - 2 e^\dagger  M^{-1} M'  M^{-1} \gamma'		\label{psi''}\\
	& \qquad - e^\dagger  M^{-1} M''  M^{-1} \gamma + 2 e^\dagger  M^{-1} M' M^{-1} M' M^{-1} \gamma,	\notag\\
u'''&= e^\dagger  M^{-1} \gamma''' - 3 e^\dagger  M^{-1} M'  M^{-1} \gamma'' - e^\dagger  M^{-1} M''' M^{-1} \gamma 		\label{psi'''}\\
&\qquad - 3 e^\dagger  M^{-1} M'' M^{-1} \gamma' +  6 e^\dagger  M^{-1} M' M^{-1} M' M^{-1} \gamma' \notag\\
&\qquad   +  3 e^\dagger  M^{-1} M'' M^{-1} M' M^{-1} \gamma + 3 e^\dagger  M^{-1} M' M^{-1} M'' M^{-1} \gamma \notag\\
&\qquad - 6 e^\dagger  M^{-1} M' M^{-1} M' M^{-1} M' M^{-1} \gamma, \notag\\
i u_s &= i e^\dagger  M^{-1} \gamma_s - i e^\dagger  M^{-1} M_s  M^{-1} \gamma, \label{psi_s} \\
u_t &= e^\dagger  M^{-1} \gamma_t - e^\dagger  M^{-1} M_t  M^{-1} \gamma. \label{psi_t}
\end{align}

Noting that $M$ is selfadjoint, we see that
\begin{align}
\overline u &= \gamma^\dagger  M^{-1} e 
\end{align}
and consequently, using \eqref{66} and the first identity in \eqref{gamma'} we get
\begin{align*}
|u|^2u &=  e^\dagger  M^{-1} \gamma \gamma^\dagger  M^{-1} e e^\dagger  M^{-1} \gamma \\
&= e^\dagger  M^{-1} [ -M'] M^{-1} [-i\Lambda M +i M \Lambda^\dagger  +  M'] M^{-1} \gamma\\
&= e^\dagger  M^{-1} M' M^{-1} \gamma'  - i e^\dagger  M^{-1} M' \Lambda^\dagger   M^{-1} \gamma - e^\dagger  M^{-1} M' M^{-1} M'  M^{-1} \gamma  .
\end{align*}

In expanding the mKdV nonlinearity, we may choose which factor of $u$ is being differentiated.  Thus, using also \eqref{66} and the first two identities in \eqref{gamma'} we obtain two different expressions: 
\begin{align*}
|u|^2u' &=   e^\dagger  M^{-1} \gamma \gamma^\dagger  M^{-1} e e^\dagger  (M^{-1} \gamma)'  \\
&= e^\dagger  M^{-1} [ -M'] M^{-1} [-i\Lambda M +i M \Lambda^\dagger  + M'] [M^{-1} \gamma' - M^{-1} M' M^{-1} \gamma ] \\
&=- e^\dagger  M^{-1} M' M^{-1} [-i\Lambda M +i M \Lambda^\dagger  + M'] M^{-1} \gamma'  \\
&\qquad + e^\dagger  M^{-1} M' M^{-1}  [-i\Lambda M +i M \Lambda^\dagger  + M'] M^{-1} M' M^{-1} \gamma \\
&=  e^\dagger  M^{-1} M' M^{-1} \gamma''
	 - i e^\dagger  M^{-1} M' \Lambda^\dagger  M^{-1} \gamma' 
	 - e^\dagger  M^{-1} M' M^{-1} M' M^{-1} \gamma'  \\
&\qquad - i e^\dagger  M^{-1} M' M^{-1} \Lambda M' M^{-1} \gamma
	 + i e^\dagger  M^{-1} M' \Lambda^\dagger   M^{-1} M' M^{-1} \gamma  \\
&\qquad	 + e^\dagger  M^{-1} M' M^{-1} M' M^{-1} M' M^{-1} \gamma,
\end{align*}
as well as
\begin{align*}
|u|^2u' &=   e^\dagger  (M^{-1} \gamma)' \gamma^\dagger  M^{-1} e e^\dagger  M^{-1} \gamma  \\
&=e^\dagger  [M^{-1} i\Lambda \gamma - M^{-1} M' M^{-1} \gamma ] \gamma^\dagger  M^{-1} [-i\Lambda M +i M \Lambda^\dagger  + M'] M^{-1} \gamma \\
&= e^\dagger  M^{-1} i\Lambda  [-M'] M^{-1} [-i\Lambda M +i M \Lambda^\dagger  + M']M^{-1} \gamma  \\
&\qquad -e^\dagger  M^{-1} M' M^{-1} [-M']  M^{-1}  [-i\Lambda M +i M \Lambda^\dagger  + M'] M^{-1} \gamma \\
&= e^\dagger  M^{-1} i\Lambda  M' M^{-1} \gamma'
	 - i e^\dagger  M^{-1} i\Lambda  M' \Lambda^\dagger  M^{-1} \gamma
	 - e^\dagger  M^{-1} i\Lambda  M' M^{-1} M' M^{-1} \gamma  \\
&\qquad - e^\dagger  M^{-1} M' M^{-1} M'  M^{-1} \gamma' 
	 + i e^\dagger  M^{-1} M' M^{-1} M' \Lambda^\dagger  M^{-1} \gamma \\
&\qquad  + e^\dagger  M^{-1} M' M^{-1} M'  M^{-1} M'  M^{-1} \gamma .
\end{align*}
Combining three copies of each expression, we deduce that
\begin{align*}
6 |u|^2u' &=  3 e^\dagger  M^{-1} M' M^{-1} \gamma''
	 + 3 e^\dagger  M^{-1} \bigl[ i\Lambda  M'  -i M'\Lambda^\dagger \bigr] M^{-1} \gamma'  \\
&\qquad	- 6e^\dagger  M^{-1} M' M^{-1} M' M^{-1} \gamma'  \\
&\qquad - 3 e^\dagger  M^{-1} M' M^{-1} \bigl[ i\Lambda M'- i M' \Lambda^\dagger \bigr] M^{-1} \gamma \\
&\qquad	- 3 e^\dagger  M^{-1} \bigl[ i\Lambda M' - i M' \Lambda^\dagger  \bigr]  M^{-1} M' M^{-1} \gamma \\
&\qquad + 3 e^\dagger  M^{-1} \Lambda  M' \Lambda^\dagger  M^{-1} \gamma
	 + 6 e^\dagger  M^{-1} M' M^{-1} M'  M^{-1} M'  M^{-1} \gamma .
\end{align*}
Recognizing the square bracket terms as $M''$ and recalling our earlier expansion \eqref{psi'''}, we deduce that
\begin{align*}
  u''' + 6|u|^2u' &= e^\dagger  M^{-1} \gamma''' - e^\dagger  M^{-1} M''' M^{-1} \gamma + 3 e^\dagger  M^{-1} \Lambda  M' \Lambda^\dagger  M^{-1} \gamma .
\end{align*}
Finally, using \eqref{M'}, \eqref{gamma_s gamma_t}, \eqref{M dot}, and \eqref{psi_t}, we conclude that 
\begin{align*}
  u_t + u''' + 6|u|^2u' &=  - e^\dagger  M^{-1} [ M_t + M''']  M^{-1} \gamma + 3 e^\dagger  M^{-1} \Lambda  M' \Lambda^\dagger  M^{-1} \gamma = 0,
\end{align*}
that is, $u$ solves \eqref{mkdv} in the variables $t,x$.

It remains to show that $u$ solves the focusing cubic nonlinear Schr\"odinger equation in the variables $s,x$.  Indeed, using  \eqref{gamma'} through \eqref{psi_s}  and the expression above for $|u|^2u$, we get
\begin{align*}
    &iu_s+u''+2|u|^2u\\
    &=i e^\dagger  M^{-1} \gamma_s - i e^\dagger  M^{-1} M_s  M^{-1} \gamma\\
    &\qquad +e^\dagger  M^{-1} \gamma'' - 2 e^\dagger  M^{-1} M'  M^{-1} \gamma'		\\
	& \qquad - e^\dagger  M^{-1} M''  M^{-1} \gamma + 2 e^\dagger  M^{-1} M' M^{-1} M' M^{-1} \gamma\\
    &\qquad+2e^\dagger  M^{-1} M' M^{-1} \gamma'  - 2i e^\dagger  M^{-1} M' \Lambda^\dagger   M^{-1} \gamma - 2e^\dagger  M^{-1} M' M^{-1} M'  M^{-1} \gamma  \\
    &=e^\dagger  M^{-1} [i\gamma_s+\gamma''] -  e^\dagger  M^{-1}[i M_s +M''] M^{-1} \gamma  - 2i e^\dagger  M^{-1} M' \Lambda^\dagger   M^{-1} \gamma   =0,
\end{align*}
which completes the proof of the proposition.
\end{proof}

\section{Asymptotics of the $N$-fold soliton}\label{S:3}

The main result in this section is the following soliton resolution statement for the solution to \eqref{mkdv} with initial data $u_0(x)=(-1)^N N \sech(x)$:

\begin{theorem}\label{T:asymp}
Fix $N\geq 1$ and let $u_N$ denote the solution to \eqref{mkdv} with initial data $u_N(x,0)=(-1)^N N\sech(x)$.  Then for each $1\leq p <\infty$, we have
\begin{align*}
\lim_{t\to \infty}\,\Bigl\|u_N(t)-(-1)^{N}\sum_{j=1}^N(2j-1)\sech\Bigl((2j-1)\big[x-(2j-1)^2t- c_j\big]\Bigr)\Bigr\|_{L^p}  = 0.
\end{align*}
The spatial shifts $c_j$ can be expressed in terms of binomial coefficients as follows:
\begin{equation}\label{c_j}
c_j=\tfrac{1}{2j-1} \ln\Bigl[\tbinom{N+j-1}{N-j} \Big/ \tbinom{2j-2}{j-1}\Bigr].
\end{equation}
\end{theorem}

Satsuma--Yajima \cite{MR463733} solved the forward scattering problem associated to any numerical multiple of the soliton, that is, for the profile $u(x)=\kappa\sech(x)$ with any value of $\kappa$.  They showed that integer multiples are exact multisolitons and they provided the associated spectral parameters.  To keep the presentation self-contained, we verify this result by direct computation.

\begin{lemma}\label{L:initial data}
Fix $N\geq 1$ and for each $1\leq j\leq N$, let $\lambda_j=i(2j-1)$ and $a_j=(-1)^j$. Let $e$, $\gamma$, and $M$ be as defined in \eqref{e, gamma} and \eqref{M} with $\theta=y=s=t=0$.  Then
\begin{align*}
  u(x)=e^\dagger  M^{-1} \gamma = (-1)^N N\sech(x).
\end{align*}
\end{lemma}

\begin{proof}
It is convenient to recognize $M$ as the Gram matrix with entries
\begin{align*}
  M_{jk} = \frac1{2} \int_{y_0}^1 y^{j-1} y^{k-1}\,dy  \qtq{with} y_0 = - e^{-2x} .
\end{align*}

Inspired by this, we introduce the upper triangular matrix $T$ with entries
\begin{equation*}
  T_{jk} = \begin{cases}(1-y_0)^{1-k} \tbinom{k-1}{j-1} (-1)^{j-1}, \qquad & \text{if}\quad  j\leq k\\0, \qquad & \text{if}\quad j>k. \end{cases}
\end{equation*}
We observe that by the binomial theorem,
\begin{align*}
  (T^\dagger  M T)_{jk} &= \frac1{2 (1-y_0)^{j+k-2}} \int_{y_0}^1  \biggl[ \sum_{m=1}^{j} \tbinom{j-1}{m-1} (-y)^{m-1} \biggr]  \biggl[ \sum_{n=1}^{k} \tbinom{k-1}{n-1} (-y)^{n-1} \biggr]\,dy  \\
  &= \frac1{2(1-y_0)^{j+k-2}} \int_{y_0}^1  (1-y)^{j+k-2} \,dy = \frac{1+e^{-2x}}{2(j+k-1)},
\end{align*}
as well as
\begin{align*}
  (T^\dagger e)_{k} &= (1-y_0)^{1-k} \sum_{j=1}^{k} \tbinom{k-1}{j-1} (-1)^{j-1} = \delta_{k1} .
\end{align*}
Noting that $\gamma_j= -y_0^{j-1}e^{-x}$, we also get
\begin{align*}
  (T^\dagger  \gamma)_{k} &=  (1-y_0)^{1-k} \sum_{j=1}^{k} \tbinom{k-1}{j-1} (-1)^{j-1} \cdot \bigl[ -y_0^{j-1}e^{-x}\bigr]= - e^{-x}.
\end{align*}

In this way, we find that 
\begin{align*}
  u(x)=e^\dagger  M^{-1} \gamma = (T^\dagger  e)^\dagger  [T^\dagger  M T]^{-1} T^\dagger  \gamma = -\frac{e^{-x}}{1+e^{-2x}}  (C^{-1} e)_1
\end{align*}
where $C$ is the Cauchy matrix with entries $C_{jk}=1/(2j+2k-2)$.  This corresponds to the choice $a_j=b_j=2j-1$ in Lemma~\ref{L:Cauchy} and therefore,
\begin{align*}
 ( C^{-1} e )_1 = 2\prod_{j=2}^N \frac{j}{1-j} = 2(-1)^{N-1} N.
\end{align*} This leads us to
\[
u(x)=(-1)^NN\frac{2e^{-x}}{1+e^{-2x}} =(-1)^NN\sech(x). \qedhere
\]
\end{proof}

The remainder of this section is dedicated to the proof of Theorem~\ref{T:asymp}. We will employ formula \eqref{E:P:psi} with the parameters given in Lemma~\ref{L:initial data}.  For ease of reference, we recapitulate this here: The solution $u_N(t,x)$ to \eqref{mkdv} with initial data
\begin{equation}\label{398}
u_N(0,x)=(-1)^N N \sech(x) \qtq{is given by} u_N(t,x)=e^\dagger M^{-1}\gamma,
\end{equation}
where $e$, $\gamma$, and $M$ have entries 
\begin{equation}\label{399}
e_j= 1, \qquad \gamma_j(t,x)=(-1)^j e^{-(2j-1)x+(2j-1)^3t}, \qtq{and} M_{jk}=\frac{1+\gamma_j\overline{\gamma}_k}{2(j+k-1)}.
\end{equation}

By contrast, the sought-after asymptotic behaviour of $u_N(t,x)$ is the sum
\begin{equation}\label{defn v_N}
v_N(t,x):=\sum_{j=1}^N \psi_j(t,x) 
\end{equation}
of single soliton solutions 
\begin{equation*}
\psi_j(t,x):=(-1)^{N}(2j-1)\sech\Bigl((2j-1)\bigl[x-(2j-1)^2t-c_j\bigr]\Bigr).
\end{equation*}
The most subtle part of Theorem~\ref{T:asymp} will be deriving the spatial shifts \eqref{c_j}; see Lemma~\ref{L:A}. 

Based on the structure of $v_N$, we partition $\R=I_0\cup I_1\cup\cdots\cup I_{2N}$, as follows:
\begin{align*}
I_m=\begin{cases}
\bigl(-\infty,\, t-\varepsilon t\bigr), &\quad\text{if }\ m=0,\\
\bigl[(2\ell-1)^2t-\varepsilon t, (2\ell-1)^2t+\varepsilon t\bigr], &\quad\text{if }\ m=2\ell-1 \text{ with $1\leq \ell\leq N$},\\
\bigl((2\ell-1)^2t+\varepsilon t, (2\ell+1)^2t-\varepsilon t\bigr), &\quad\text{if }\ m=2\ell   \text{ with $1\leq \ell\leq N-1$},\\
\bigl((2N-1)^2t+\varepsilon t,\infty\bigr), &\quad\text{if }\ m=2N,
\end{cases}
\end{align*}
for some fixed $\varepsilon\in (0,\tfrac{1}{2N}]$.  Throughout this section, the constants implicit in the notation  $\lesssim$ will be allowed to depend on $N$ and the choice of $\varepsilon$.

The proof of Theorem~\ref{T:asymp} will be comprised of two parts: Lemma~\ref{L:even intervals}, which handles the contribution from the even-numbered intervals, and Lemma~\ref{L:odd intervals}, which covers the odd-numbered intervals.  

Each odd-numbered interval $I_{2\ell-1}$ follows the corresponding soliton $\psi_\ell(t,x)$.  Their widths grow linearly in time, so we need not worry that they are not precisely centered (due to the shifts $c_j$).  The analysis of these intervals relies on subtle cancellations between $u_N$ and $v_N$.

The even numbered intervals $I_{2\ell}$, $0\leq \ell\leq N$, occupy the space between the solitons, where both $u_N$ and $v_N$ are asymptotically negligible.  These are easier to treat:


\begin{lemma}\label{L:even intervals}
Fix $1\leq p< \infty$. For each $0\leq \ell\leq N$, we have
\begin{align}\label{goal even}
\lim_{t\to \infty}\, \bigl\|1_{I_{2\ell}}(x)\bigl[u_N(t,x)-v_N(t,x)\bigr]\bigr\|_{L^p_x}=0.
\end{align}
\end{lemma}

\begin{proof} For $1\leq j\leq N$, a simple change of variables shows that
\begin{align*}
\bigl\| 1_{\R\setminus I_{2j-1}(t)}(x)\psi_j(t,x)\bigr\|_{L^p_x} = (2j-1)^{1-\frac1p}\bigl\| 1_{\R\setminus [-(2j-1)\varepsilon t, (2j-1)\varepsilon t]}(x)\sech(x-c_j)\bigr\|_{L^p_x}
\end{align*}
and consequently, 
\begin{align}\label{4:56}
\bigl\| 1_{\R\setminus I_{2j-1}(t)}(x)\psi_j(t,x)\bigr\|_{L^p_x} \lesssim e^{-\varepsilon t} .
\end{align}

Using the triangle inequality and \eqref{4:56}, we immediately deduce that
$$
\lim_{t\to \infty}\, \bigl\|1_{I_{2\ell}}(x)\, v_N(t,x)\bigr\|_{L^p_x}=0.
$$
Thus \eqref{goal even} will follow once we demonstrate
\begin{align}\label{goal even'}
\lim_{t\to \infty}\, \bigl\|1_{I_{2\ell}}(x) \,u_N(t,x)\bigr\|_{L^p_x}=0.
\end{align}

Working towards this goal, we note that whenever $x\in I_{2\ell}$ with $0\leq \ell\leq N$ and $t\geq 0$, we have
\begin{align}\label{asymp gamma 2J}
\begin{cases}
|\gamma_j(t,x)|\leq e^{-x+(2\ell-1)^2t}\leq e^{-\varepsilon t},  & \quad\text{for } \ j\leq \ell, \\
|\gamma_j(t,x)|\geq e^{-x+(2\ell+1)^2 t}\geq e^{\varepsilon t} ,  & \quad \text{for }\ j\geq \ell+1.
\end{cases}
\end{align}
Inspired by this, we write
\begin{align}\label{rewrite u}
u_N=e^\dagger M^{-1}\gamma
&=e^\dagger  \dg{1,\ldots,1,  \gamma_{\ell+1}^{-1},\ldots, \gamma_{N}^{-1} }\dg{1,\ldots,1,  \gamma_{\ell+1},\ldots, \gamma_{N} }M^{-1}\notag\\
&\qquad    \dg{1,\ldots,1,\gamma_{\ell+1}, \ldots, \gamma_N}\dg{1,\ldots,1,\gamma_{\ell+1}^{-1}, \ldots, \gamma_N^{-1}}\gamma  \\
&=(1,\ldots, 1, \gamma_{\ell+1}^{-1},\ldots, \gamma_{N}^{-1}) W^{-1} (\gamma_1, \ldots, \gamma_{\ell}, 1,\ldots, 1)^\dagger, \notag
\end{align}    
where $W$ is the $N\times N$ matrix
$$
W= \dg{1,\ldots,1,  \gamma_{\ell+1}^{-1},\ldots, \gamma_{N}^{-1} }\, M\, \dg{1,\ldots,1,  \gamma_{\ell+1}^{-1},\ldots, \gamma_{N}^{-1} }.
$$

A simple computation reveals that the entries of the $N\times N$ matrix $W$ are as follows:
\begin{equation}\label{W2}\newcommand{\tmp}[1]{\vphantom{\biggl|}\tfrac{1}2\cdot\tfrac{#1}{j+k-1}}\begin{tabular}{|c|c|c|c|}
\hline
$W_{jk}$ & $1\leq k \leq\ell$ & $\ell< k \leq N$ \\
\hline
$1\leq j \leq \ell$ &  $\tmp{1+\gamma_j\gamma_k}$ & $\tmp{\gamma_k^{-1}+\gamma_j}$ \\
\hline
$\ell< j \leq N$  &  $\tmp{\gamma_j^{-1} + \gamma_k}$ & $\tmp{1+\gamma_j^{-1}\gamma_k^{-1}}$ \\
\hline
\end{tabular}
\end{equation}
In view of \eqref{asymp gamma 2J}, the $t\to\infty$ limit of $W$ has the form
\begin{align*}
A=\begin{pmatrix}
    A_1 & \bf0\\
    \bf0 & A_2 \end{pmatrix} \ \ \text{with}\ \  A_1=\Bigl(\tfrac 1{2(j+k-1)}\Bigr)_{ 1\leq j,k\leq \ell } \ \ \text{and}\ \  A_2=\Bigl(\tfrac 1{2(j+k-1)}\Bigr)_{\ell+1\leq j,k\leq N}.
\end{align*}
Moreover, the remainder $O$ in the decomposition $W=A+O$ satisfies 
\begin{align}\label{W decomp}
\|O\|\lesssim \sup_{1\leq j\leq\ell} \bigl| \gamma_j\bigr| + \sup_{\ell<k<N} \bigl| \gamma_k^{-1}\bigr| \lesssim e^{-\varepsilon t}.
\end{align}

Both $A_1$ and $A_2$ are Cauchy matrices.  By Lemma~\ref{L:Cauchy}, they are invertible.  Thus, the matrix $A$ is also invertible.
Using \eqref{W decomp}, we conclude that there exists $t_0=t_0(\varepsilon)$ sufficiently large so that the matrix $W$ is invertible for $t\geq t_0$ and
\begin{align}\label{W norm}
\|W^{-1}\|\leq 2\|A^{-1}\| \lesssim 1 \qtq{for} t\geq t_0.
\end{align}    
Using this, \eqref{W decomp}, and the resolvent identity, we obtain
\begin{align*}
\|W^{-1} -A^{-1}\| = \|A^{-1} O W^{-1}\|\lesssim e^{-\varepsilon t} \qtq{for} t\geq t_0.
\end{align*}
Recalling \eqref{rewrite u}, for $x\in I_{2\ell}$ we may thus bound
\begin{align}\label{uN bdd}
|u_N(t,x)|&= \bigl|(1,\ldots, 1,\  \gamma_{\ell+1}^{-1},\ldots, \gamma_{N}^{-1}) W^{-1} (\gamma_1, \ldots, \gamma_{\ell}, 1,\ldots, 1)^{\dagger}\bigr|\notag\\
&\lesssim  \bigl|(1,\ldots, 1,\  \gamma_{\ell+1}^{-1},\ldots, \gamma_{N}^{-1}) A^{-1} (\gamma_1, \ldots, \gamma_{\ell}, 1,\ldots, 1)^{\dagger}\bigr|
	+ e^{-\varepsilon t}\notag\\
&\lesssim \bigl|(1,\ldots, 1) A_1^{-1} (\gamma_1, \ldots, \gamma_{\ell})^{\dagger} \bigr|
	+ \bigl|(\gamma_{\ell+1}^{-1},\ldots, \gamma_{N}^{-1}) A_2^{-1} (1,\ldots, 1)^{\dagger}\bigr|+ e^{-\varepsilon t}\\
&\lesssim  e^{-\varepsilon t}, \notag
\end{align}
where in the last step we used \eqref{asymp gamma 2J} once more.  

When $1\leq \ell\leq N-1$, employing \eqref{uN bdd} we find
\begin{align*}
\bigl\|1_{I_{2\ell}}(x) \,u_N(t,x)\bigr\|_{L^p_x} \lesssim e^{-\varepsilon t} |I_{2\ell}|^{\frac1p} \lesssim t^{\frac1p} e^{-\varepsilon t} \to 0 \qtq{as} t\to \infty,
\end{align*}
which settles \eqref{goal even'} in this case.

When $\ell=0$, we use \eqref{asymp gamma 2J}, \eqref{rewrite u}, and \eqref{W norm} to obtain
\begin{align*}
|u_N(t,x)|&=\bigl|(\gamma_1^{-1},\ldots, \gamma_N^{-1})W^{-1}(1,\ldots,1)^\dagger\bigr|\lesssim \|W^{-1}\|\sup_{1\leq j\leq N} |\gamma_j(t,x)^{-1}|\lesssim e^{x-t} 
\end{align*}
for $x\in I_0$ and $t\geq t_0$.  Thus,
\begin{align*}
\bigl\|1_{I_0}(x) \,u_N(t,x)\bigr\|_{L^p_x} \lesssim_p e^{t - \varepsilon t} e^{-t} \to 0 \qtq{as} t\to \infty,
\end{align*}
yielding \eqref{goal even'} for $\ell=0$.

Finally, for $\ell=N$ we use again \eqref{asymp gamma 2J}, \eqref{rewrite u}, and \eqref{W norm} to obtain
\begin{align*}
|u_N(t,x)|&=\bigl|(1,\ldots,1)W^{-1}(\gamma_1,\ldots, \gamma_N)^\dagger\bigr|\lesssim \|W^{-1}\| \sup_{1\leq j\leq N} |\gamma_j(t,x)| \lesssim e^{-x+(2N-1)^2t}
\end{align*}
for $x\in I_{2N}$ and $t\geq t_0$. This yields
\begin{align*}
\bigl\|1_{I_{2N}}(x) \,u_N(t,x)\bigr\|_{L^p_x} \lesssim_p e^{-(2N-1)^2t - \varepsilon t} e^{(2N-1)^2t} \to 0 \qtq{as} t\to \infty,
\end{align*}
which proves \eqref{goal even'} for $\ell=N$.
\end{proof}

In the previous proof, we saw that after suitable conjugation, the product $e^\dagger M^{-1}\gamma$ could be written as a main term, which involved the matrix $A$ and an error that is exponentially smaller.  A similar argument prevails in the case of the intervals $I_{2\ell-1}$. However, the evaluation of the main term is much more complicated; this is the subject of Lemma~\ref{L:A}.

\begin{lemma}\label{L:odd intervals}
Fix $1\leq p< \infty$. For each $1\leq \ell\leq N$, we have
\begin{align}\label{goal odd}
\lim_{t\to \infty}\, \bigl\|1_{I_{2\ell-1}}(x)\bigl[u_N(t,x)-v_N(t,x)\bigr]\bigr\|_{L^p_x}=0.
\end{align}
\end{lemma}

\begin{proof}
Fix $1\leq \ell\leq N$.  By \eqref{4:56} we have
$$
\lim_{t\to\infty}\,\bigl\|1_{I_{2\ell-1}}(x)\psi_j(t,x)\bigr\|_{L^p_x}=0 \qtq{whenever} j\neq \ell.
$$
Thus, by the triangle inequality, it suffices to show that 
\begin{align}\label{goal odd'}
\lim_{t\to \infty}\, \bigl\|1_{I_{2\ell-1}}(x)\bigl[u_N(t,x)-\psi_\ell(t,x)\bigr]\bigr\|_{L^p_x}=0.
\end{align}

For $x\in I_{2\ell-1}$ and $t\geq 0$, we have 
\begin{align}\label{asymp gamma 2J-1}
\begin{cases}
|\gamma_j(t,x)|\leq e^{-x+(2\ell-3)^2t}\leq e^{-7t},   &\quad \text{for }\ j<\ell, \\
|\gamma_j(t,x)|\geq e^{-x+(2\ell+1)^2t}\geq e^{7t},  &\quad \text{for }\ j>\ell,\\
  e^{-t}\leq  e^{-(2\ell-1)\varepsilon t}\leq |\gamma_\ell(t,x)|\leq  e^{(2\ell-1)\varepsilon t}\leq e^{t}. &  
\end{cases}
\end{align}
As a consequence, we also have 
\begin{align}\label{asymp gamma 2J-1'}
\begin{cases}
|\gamma_j(t,x)\gamma_\ell(t,x)|\leq e^{-6t},  &\quad \text{for }\ j<\ell,  \\
|\gamma_j(t,x)\gamma_\ell(t,x)|\geq e^{6t}, &\quad \text{for }\ j>\ell. 
\end{cases}
\end{align}

Inspired by this, we write
\begin{align}\label{rewrite u odd}
u_N= e^\dagger M^{-1}\gamma
&=e^\dagger  \dg{1,\ldots,1,  \gamma_{\ell+1}^{-1},\ldots, \gamma_{N}^{-1} }\dg{1,\ldots,1,  \gamma_{\ell+1},\ldots, \gamma_{N} }M^{-1}\notag\\
&\qquad \dg{1,\ldots,1, \gamma_{\ell}, \ldots, \gamma_N}\dg{1,\ldots,1, \gamma_{\ell}^{-1}, \ldots, \gamma_N^{-1}}\gamma\notag\\
&=(1,\ldots, 1,  \gamma_{\ell+1}^{-1},\ldots, \gamma_{N}^{-1}) W^{-1} (\gamma_1, \ldots, \gamma_{\ell-1}, 1,\ldots, 1)^\dagger,
\end{align}
with 
\begin{align*}
W= \dg{1,\ldots,1,  \gamma_{\ell}^{-1},\ldots, \gamma_{N}^{-1} }\, M\, \dg{1,\ldots,1,  \gamma_{\ell+1}^{-1},\ldots, \gamma_{N}^{-1} }.
\end{align*}

Direct computation reveals that the entries of the $N\times N$ matrix $W$ are as follows:
\begin{equation}\label{W3}\newcommand{\tmp}[1]{\vphantom{\biggl|}\tfrac{1}2\cdot\tfrac{#1}{j+k-1}}\begin{tabular}{|c|c|c|c|}
\hline
$W_{jk}$ & $1\leq k < \ell$ & $k =\ell$ & $\ell< k \leq N$ \\
\hline
$1\leq j < \ell$ & $\tmp{1+\gamma_j\gamma_k}$ & $\tmp{1+\gamma_j\gamma_k}$ & $\tmp{\gamma_k^{-1}+\gamma_j}$ \\
\hline
$j=\ell$            & $\tmp{\gamma_j^{-1}+\gamma_k}$ & $\tmp{\gamma_j^{-1}+ \gamma_k}$ & $\tmp{1+\gamma_j^{-1}\gamma_k^{-1}}$ \\
\hline
$\ell< j \leq N$  & $\tmp{\gamma_j^{-1}+\gamma_k}$ & $\tmp{\gamma_j^{-1} + \gamma_k}$ & $\tmp{1+\gamma_j^{-1}\gamma_k^{-1}}$ \\
\hline
\end{tabular}
\end{equation}
In view of \eqref{asymp gamma 2J-1} and \eqref{asymp gamma 2J-1'}, we may thus decompose 
\begin{align}\label{W decomp odd}
W=A+O \qtq{with}  \|A\|\lesssim e^{t} \qtq{and} \|O\|\lesssim e^{-6t},
\end{align}
where the $N\times N$ matrix $A$ has the following form:
\begin{equation}\label{evil A}
\newcommand{\tmp}[1]{\vphantom{\biggl|}\tfrac{1}2\cdot\tfrac{#1}{j+k-1}}\begin{tabular}{|c|c|c|c|}
\hline
$A_{jk}$ & $1\leq k < \ell$ & $k =\ell$ & $\ell< k \leq N$ \\
\hline
$1\leq j < \ell$ & $\tmp{1}$ & $\tmp{1}$ & $0$ \\
\hline
$j=\ell$            & $\tmp{\gamma_\ell^{-1}}$ & $\tmp{\gamma_\ell^{-1}+ \gamma_\ell^{ }}$ & $\tmp{1}$ \\
\hline
$\ell< j \leq N$  & $0$ & $\tmp{ \gamma_\ell}$ & $\tmp{1}$ \\
\hline
\end{tabular}
\end{equation}

\begin{lemma}\label{L:A}
Fix $1\leq \ell\leq N$ and let $A$ denote the matrix with entries \eqref{evil A}.
Then $A$ is invertible and
\begin{align}\label{A-1 bdd}
\|A^{-1}\|\lesssim 1 \qtq{uniformly for} t,x\in\R.
\end{align}
Furthermore, the shifted soliton $\psi_\ell(t,x)$ admits the representation 
\begin{align}\label{A psi}
\psi_\ell = (\underbrace{1, \ldots, 1}_{\ell},  \underbrace{0,\dots, 0}_{N-\ell})
	A^{-1} (\underbrace{0, \ldots, 0}_{\ell-1},\underbrace{1,\ldots, 1}_{N+1-\ell})^\dagger.
\end{align}
\end{lemma}

Postponing the proof of this lemma for now, let us see how it can be used to complete the proof of Lemma~\ref{L:odd intervals}.

Combining Lemma~\ref{L:A} with \eqref{W decomp odd}, we deduce that there exists $t_1=t_1(\varepsilon)$ so that for $t\geq t_1$ the matrix $W$ is invertible and
$$
\|W^{-1}\|\leq 2\|A^{-1}\|\lesssim 1.
$$
Using this, the resolvent identity, and \eqref{W decomp odd}, we obtain
$$
\|W^{-1}-A^{-1}\|= \|A^{-1} O W^{-1}\|\lesssim e^{-6t} \qtq{for} t\geq t_1.
$$
Together with \eqref{asymp gamma 2J-1}, this yields
\begin{align}\label{error}
\bigl|(1,\ldots 1,  \gamma_{\ell+1}^{-1},\ldots, \gamma_{N}^{-1})( W^{-1}-A^{-1}) (\gamma_1, \ldots, \gamma_{\ell-1}, 1,\ldots, 1)^\dagger\bigr|\lesssim e^{-6t}
\end{align}
whenever $t\geq t_1$ and $x\in I_{2\ell-1}$.

On the other hand, using Lemma~\ref{L:A} and \eqref{asymp gamma 2J-1} we get
\begin{align}\label{main}
\bigl|(1,\ldots, 1,  \gamma_{\ell+1}^{-1},\ldots, &\gamma_{N}^{-1})A^{-1} (\gamma_1, \ldots, \gamma_{\ell-1}, 1,\ldots 1)^\dagger-\psi_\ell\bigr|\notag\\
&\lesssim \|A^{-1}\| \Bigl[ \, \sup_{\ell+1\leq j\leq N}|\gamma_j^{-1}| + \sup_{1\leq j\leq \ell-1} |\gamma_j|\Bigr]\lesssim e^{-7t}
\end{align}
for all $t\geq t_1$ and $x\in I_{2\ell-1}$.

Combining \eqref{rewrite u odd} with \eqref{error} and \eqref{main}, we deduce that 
\begin{align*}
\bigl\| 1_{I_{2\ell-1}}(x)\bigl[u_N(t,x) - \psi_\ell(t,x)\bigr]\bigr\|_{L^p_x}\lesssim |I_{2\ell-1}|^{\frac1p} e^{-6t} \lesssim t^{\frac1p} e^{-6t} \to 0 \qtq{as} t\to \infty,
\end{align*}
which settles \eqref{goal odd'} and completes the proof of the lemma.
\end{proof}

To complete the proof of Theorem~\ref{T:asymp}, it remains to prove Lemma~\ref{L:A}.

\begin{proof}[Proof of Lemma~\ref{L:A}]
We write the matrix $A$ as
$$A=\begin{pmatrix}
      B & \,  b  & \, \bf0\\[2mm]
    \gamma_\ell^{-1}\, b^\dagger  &\, (\gamma_\ell^{-1}+\gamma_\ell) d & \, f^\dagger \\[2mm]
     \bf0 & \gamma_\ell \,f &\, H
 \end{pmatrix},$$
where 
\begin{align*}
B= \Bigl(\tfrac{1}{2(j+k-1)} \Bigr)_{1\leq j,k<\ell}, &\qquad H = \Bigl(\tfrac{1}{2(j+k-1)} \Bigr)_{\ell< j,k\leq N},\\
b= \Bigl(\tfrac{1}{2(j+\ell-1)}\Bigr)_{1\leq j<\ell}, \qquad f= &\Bigl(\tfrac{1}{2(\ell+j-1)}\Bigr)_{\ell<j\leq N}, \qquad d= \tfrac{1}{2(\ell+\ell-1)}.
\end{align*}

Performing Gaussian elimination, as in the derivation of the Schur complement, we obtain the following formula for the inverse of $A$:
\begin{align}\label{A-1}
A^{-1}=\begin{pmatrix}
B^{-1} + z\gamma_\ell^{-1}B^{-1}bb^\dagger B^{-1} & -zB^{-1}b &  zB^{-1}bf^\dagger H^{-1}\\[2mm]
     -z\gamma_\ell^{-1} b^\dagger B^{-1} & z & -zf^\dagger H^{-1}\\[2mm]
     z  H^{-1}fb^\dagger B^{-1} & -z\gamma_\ell H^{-1}f & H^{-1}+ z\gamma_\ell H^{-1}ff^\dagger H^{-1}
 \end{pmatrix}
 \end{align}
where
\begin{equation}\label{z defn}
z:= \Bigl\{\gamma_\ell^{-1}\bigl[d- b^\dagger  B^{-1}b\bigr] + \gamma_\ell\bigl[d-  f^\dagger H^{-1} f\bigr]\Bigr\}^{-1} .
\end{equation}
In order to make this rigorous, we must confirm that $B$ and $H$ are actually invertible and likewise, that the definition of $z$ does not involve division by zero.

As $B$ and $H$ are Cauchy matrices, Lemma~\ref{L:Cauchy} guarantees that they are invertible. We now turn our attention to $z$.  

By Lemma~\ref{L:Cauchy}, we know that the Cauchy matrices
\begin{align*}
C_1:=\begin{pmatrix} B & b\\[2mm]  b^\dagger  & d  \end{pmatrix} = \Bigl(\tfrac{1}{2(j+k-1)}\Bigr)_{1\leq j,k\leq\ell}
\end{align*}
and 
\begin{align*}
C_2:=\begin{pmatrix} d & f^\dagger \\[2mm]  f& H  \end{pmatrix} =\Bigl(\tfrac{1}{2(j+k-1)}\Bigr)_{\ell\leq j,k\leq N}.
\end{align*}
are both invertible.  Thus, we may apply the Schur complement to obtain
\begin{align}\label{8:44}
(d- b^\dagger  B^{-1}b)^{-1}&=(0,\ldots, 0 , 1)\begin{pmatrix}
 B & b\\ b^\dagger  & d
\end{pmatrix}^{-1}(0,\ldots, 0 , 1)^\dagger \notag\\
&=\,(C_1^{-1})_{\ell\ell}= 2(2\ell-1)\tbinom{2\ell-2}{\ell-1}^2
\end{align}
and 
\begin{align}\label{8:45}
(d-  f^\dagger H^{-1} f)^{-1}&= (1, 0,\ldots, 0 )\begin{pmatrix}
 d & f^\dagger \\ f & H
\end{pmatrix}^{-1}(1, 0,\ldots, 0)^\dagger \notag\\
&= (C_2^{-1})_{\ell \ell} = 2(2\ell-1)\tbinom{N+\ell-1}{2\ell-1}^2 , 
\end{align}
where the last equalities in \eqref{8:44} and \eqref{8:45} follow from Lemma~\ref{L:Cauchy}.  Crucially, these two values are both positive and consequently,
$$
\Bigl| \gamma_\ell^{-1}\bigl[d- b^\dagger  B^{-1}b\bigr] + \gamma_\ell\bigl[d-  f^\dagger H^{-1} f\bigr] \Bigr|
	\simeq |\gamma_\ell|^{-1} + |\gamma_\ell| \gtrsim 1.
$$
This allows us not only to define $z$, but also provides the bounds
\begin{equation}\label{z}
 |z| + |z \gamma_\ell|  + |z \gamma_\ell^{-1}| \lesssim 1 \qtq{uniformly for} t,x \in\R.
\end{equation}
In view of \eqref{A-1}, this proves \eqref{A-1 bdd}.

It remains to prove \eqref{A psi}.  From \eqref{A-1} and elementary manipulations, we find
 \begin{align}\label{8:41}
(\underbrace{1,\ldots, 1}_{\ell\text{ times}} &, 0,\dots, 0) A^{-1} (\underbrace{0, \ldots, 0}_{\ell-1\text{ times}},1,\dots, 1)^\dagger\notag\\
&\quad=(1,\ldots,1)
       \begin{pmatrix}
           -zB^{-1}b & zB^{-1}bf^\dagger H^{-1}\\
           z & -zf^\dagger H^{-1}
           \end{pmatrix}
       (1,\ldots, 1)^\dagger \notag\\
&\quad=z(1-e^\dagger B^{-1}b-f^\dagger H^{-1}e+e^\dagger B^{-1}bf^\dagger H^{-1}e )\notag\\
&\quad=z(1-e^\dagger B^{-1}b)(1-f^\dagger H^{-1}e)\notag\\
&\quad=\tfrac{1-e^\dagger B^{-1}b}{d-b^\dagger B^{-1}b}\cdot\tfrac{1-f^\dagger  H^{-1}e}{d- f^\dagger H^{-1}f}
	\Bigl[\gamma_\ell^{-1}(d-f^\dagger H^{-1}f)^{-1}+\gamma_\ell(d-b^\dagger B^{-1}b)^{-1}\Bigr]^{-1}. 
\end{align}
The value of the final factor in \eqref{8:41} is apparent from \eqref{8:44} and \eqref{8:45}.  The values of the first two factors can be found by applying Lemma~\ref{L:Cauchy} and the Schur complement to the matrices $C_1$ and $C_2$ introduced earlier:\begin{align}\label{8:42}
\frac{1-e^\dagger B^{-1}b}{d-b^\dagger B^{-1}b}
&=e^\dagger  \begin{pmatrix}
     B & b\\ b^\dagger & d
 \end{pmatrix}^{-1}(0,\dots, 0,  1)^\dagger \notag\\
 &= e^\dagger  C_1^{-1}(0,\dots, 0, 1)^\dagger = 2(2\ell-1) \tbinom{2\ell-2}{\ell-1}
\end{align}
and 
\begin{align}\label{8:43}
\frac{1-f^\dagger H^{-1}e}{d- f^\dagger H^{-1}f}
 &=(1, 0,\dots, 0)\begin{pmatrix}
     d & f^\dagger \\  f& H
 \end{pmatrix}^{-1}e\notag\\
 &=(1, 0,\ldots, 0)C_2^{-1}e= 2 (-1)^{N-\ell} (2\ell-1) \tbinom{N+\ell-1}{2\ell-1}.  
\end{align}
Notice the similarity between these expressions and those in \eqref{8:44} and \eqref{8:45}, which we use to evaluate the last factor in \eqref{8:41}.  
This leads to a dramatic simplification:
\begin{align*}
(\underbrace{1, \ldots, 1}_{\ell},  \underbrace{0,\dots, 0}_{N-\ell}) A^{-1} (\underbrace{0, \ldots, 0}_{\ell-1},\underbrace{1,\dots, 1}_{N+1-\ell})^\dagger
= (-1)^{N-\ell} (2\ell-1) \frac{2}{\widetilde a_\ell\gamma_\ell+(\widetilde a_\ell\gamma_\ell)^{-1}},
\end{align*}
where
$$
\widetilde a_\ell=\tbinom{2\ell-2}{\ell-1}\Big/ \tbinom{N+\ell-1}{N-\ell}.
$$
This is simply a rewriting of our claim \eqref{A psi}.
\end{proof}

We end this section with the following corollary of Theorem~\ref{T:asymp}, which proves the soliton resolution statement in Fourier--Lebesgue spaces:

\begin{corollary}\label{C:asymp}
Fix $N\geq 1$ and let $u_N$ denote the solution to \eqref{mkdv} with initial data $u_N(x,0)=(-1)^N N\sech(x)$.  Then for each $1\leq p <\infty$ and $s\in \R$, we have
\begin{align*}
\lim_{t\to \infty}\,\Bigl\|u_N(t)-(-1)^{N}\sum_{j=1}^N(2j-1)\sech\Bigl((2j-1)\big[x-(2j-1)^2t- c_j\big]\Bigr)\Bigr\|_{\mathcal{F}L^p_s(\R)}  = 0,
\end{align*}
where the spatial shifts $c_j$ are given by \eqref{c_j}.
\end{corollary}

\begin{proof}
As in the proof of Theorem~\ref{T:asymp}, we write
\begin{equation*}
v_N(t,x)=\sum_{j=1}^N \psi_j(t,x) 
\end{equation*}
with 
\begin{equation*}
\psi_j(t,x)=(-1)^{N}(2j-1)\sech\Bigl((2j-1)\bigl[x-(2j-1)^2t-c_j\bigr]\Bigr).
\end{equation*}
Thus, the claim reduces to proving 
\begin{align}\label{12:58}
\lim_{t\to \infty}\,\bigl\|u_N(t)-v_N(t)\bigr\|_{\mathcal{F}L^p_s(\R)}  = 0 \qtq{for all} 1\leq p<\infty \qtq{and} s\in \R.
\end{align}
Clearly, \eqref{12:58} for one value of $s$ implies this statement for all smaller values of $s$; thus, it suffices to prove \eqref{12:58} for $s\in \mathbb{N}$.

Recall now the classic inductive argument of Lax \cite{MR369963}.  As noted in \eqref{Hs int}, this guarantees that
\begin{align}\label{12:59}
\sup_{t\in \R}\, \bigl\|u_N(t)\bigr\|_{\mathcal{F}L^2_m(\R)}\leq C\bigl(  \bigl\|u_N(0)\bigr\|_{\mathcal{F}L^2_m(\R)}\bigr)
	\lesssim_{N,m} 1 \qtq{for integers} m\geq 0.
\end{align}

A simple computation yields
\begin{align}\label{12:60}
\bigl\|v_N(t)\bigr\|_{\mathcal{F}L^2_m(\R)}\leq \sum_{j=1}^N \bigl\|\psi_j(t)\bigr\|_{\mathcal{F}L^2_m(\R)}\lesssim_m\sum_{j=1}^N (2j-1)^{m+\frac12} \lesssim_m N^{m+\frac32},
\end{align}
uniformly for $t\in\R$.

Combining \eqref{12:59} and \eqref{12:60} we conclude
\begin{align}\label{12:61}
\sup_{t\in \R}\, \bigl\|u_N(t)\bigr\|_{\mathcal{F}L^2_m(\R)} +\sup_{t\in \R}\, \bigl\|v_N(t)\bigr\|_{\mathcal{F}L^2_m(\R)}\lesssim_{N,m}1 \quad\text{for all integers $m\geq 0$.}
\end{align}

First, we establish \eqref{12:58} for $p=1$ and $s\in \mathbb{N}$. By Cauchy--Schwarz and \eqref{12:61}, 
\begin{align*}
\bigl\|u_N(t)&-v_N(t)\bigr\|_{\mathcal{F}L^1_s(\R)}\\
&\leq \bigl\|\langle \xi\rangle^s\bigl[\widehat{u_N}(t)-\widehat{v_N}(t)\bigr]\bigr\|_{L^1(|\xi|\le M)}+\bigl\|\langle \xi\rangle^s\bigl[\widehat{u_N}(t)-\widehat{v_N}(t)\bigr]\bigr\|_{L^1(|\xi|\ge M)}\\
&\lesssim  M^{s+\frac12} \bigl\|u_N(t)-v_N(t)\bigr\|_{L^2} + \bigl\|u_N(t)-v_N(t)\bigr\|_{\mathcal{F}L^2_{s+1}(\R)} \bigl\|\langle \xi\rangle^{-1}\bigr\|_{L^2(|\xi|\geq M)}\\
&\lesssim M^{s+\frac12}\bigl\|u_N(t)-v_N(t)\bigr\|_{L^2}+M^{-\frac{1}{2}} \Bigl[\bigl\|u_N(t)\bigr\|_{\mathcal{F}L^2_{s+1}(\R)}+ \bigl\|v_N(t)\bigr\|_{\mathcal{F}L^2_{s+1}(\R)}\Bigr]\\
&\lesssim_{N,s} M^{s+\frac12}\bigl\|u_N(t)-v_N(t)\bigr\|_{L^2}+M^{-\frac{1}{2}}.
\end{align*}
By Theorem~\ref{T:asymp}, this converges to zero as $t, M\to \infty$, yielding \eqref{12:58} for $p=1$.

Finally, for $1<p<\infty$ and $s\in \mathbb{N}$, we may simply bound
\begin{align*}
\bigl\|u_N(t)-v_N(t)\bigr\|_{\mathcal F L_s^p(\R)}
&\le\bigl\|\widehat{u_N}(t)-\widehat{v_N}(t)\bigr\|_{L^\infty}^{\frac{p-1}p} \bigl\|\langle \xi\rangle^{ps}\bigl[\widehat{u_N}(t)-\widehat{v_N}(t)\bigr]\bigr\|_{L^1}^{\frac1p}\\
&\lesssim \bigl\|u_N(t)-v_N(t)\bigr\|_{L^1}^{\frac{p-1}p} \bigl\|u_N(t)-v_N(t)\bigr\|_{\mathcal{F}L^1_{\ceil{ps}}(\R)}^{\frac1p}.
\end{align*}
This converges to zero as $t\to \infty$ by Theorem~\ref{T:asymp} and the case $p=1$ of \eqref{12:58}.
\end{proof}

\section{Unbounded growth of Fourier--Lebesgue norms}\label{S:4}

This section is dedicated to the proof of Theorem~\ref{T:main}.  The key ingredient is Proposition~\ref{P:uN}, which provides bounds on Fourier--Lebesgue norms of solutions to \eqref{mkdv} with initial data that represents spatial rescalings of $u_0(x)= (-1)^N N\sech(x)$.  We will derive Theorem~\ref{T:main} from Proposition~\ref{P:uN} using the time-reversibility and time-translation symmetries of \eqref{mkdv}. 

\begin{proposition}\label{P:uN}
Given  $N\in\nN$ and $\lambda \geq N$, let $u_{N,\lambda}$ denote the solution to \eqref{mkdv} with initial data $u_{N, \lambda}(0,x)=(-1)^N N\lambda^{-1}\sech\bigl(\frac x\lambda\bigr)$.
For each $1\leq  p<\infty$ and $s\in \R$, we have
\begin{align}\label{uN0}
\bigl\|u_{N,\lambda}(0)\bigr\|_{\mathcal{F}L^p_s(\R)}\simeq_{p,s} N \lambda^{-\frac1p}.
\end{align}
Moreover, if $p\neq2$, then there is an $\alpha(p)>0$ and a choice of $t_*>0$, which may depend on $p$, $s$, $N$, and $\lambda$, so that
\begin{alignat}{2}
\bigl\|u_{N,\lambda}(t_*)\bigr\|_{\mathcal{F}L^p_s(\R)} &\lesssim_{p,s} N^{1-\alpha(p)} \lambda^{-\frac1p} &\qtq{when}  &2<p<\infty,\text{ and} \label{uN small} \\
\bigl\|u_{N,\lambda}(t_*)\bigr\|_{\mathcal{F}L^p_s(\R)} &\gtrsim_{p,s} N^{1+\alpha(p)} \lambda^{-\frac1p} &\qtq{when}  &1\leq p<2.  \label{uN big}
\end{alignat}
Importantly, the implicit constant is independent of $\lambda$ and $N$.
\end{proposition}

The time $t_*$ will be chosen to be so large that $u_{N,\lambda}(t_*)$ has resolved into well separated solitons, as described in Theorem~\ref{T:asymp}.  Taking the Fourier transform, we see that $\widehat{u_{N,\lambda}}(t_*)$ is comprised of $N$ wavepackets with widely separated characteristic frequencies.  We will control the resulting sum using the following result of Bourgain~\cite{MR1029904}.

\begin{lemma}[\cite{MR1029904}]\label{L:Bourgain}
For any $N\in \nN$, $a= (a_1, a_2, \ldots, a_N)\in \C^N$, and $2<p<\infty$,
\begin{equation}\label{E:Bourgain}
\Bigl\|\sum_{n=1}^N a_n e^{in^2 t}\Bigr\|_{L^p_t([0,2\pi])}^p\lesssim_p N^{\theta(p)}\| a\|_{\l^2}^p
\end{equation}
with $\theta(p)= \frac p2- 2$ when $p>4$ and with any  $0<\theta(p)<\tfrac p2 -1$ when $2<p\leq4$.
\end{lemma}

\begin{proof}
For $p>4$, this was proved in \cite[Proposition~1.10]{MR1029904}.  This paper also shows that there exists a universal $c>0$ so that
$$
\Bigl\|\sum_{n=1}^N a_n e^{in^2 t}\Bigr\|_{L^4_t([0,2\pi])}\lesssim e^{c\frac{\log N}{\log\log N}}\| a\|_{\l^2},
$$
uniformly for $N\geq1$ and $a\in \C^N$.  

On the other hand, it is trivial to see that 
$$
\Bigl\|\sum_{n=1}^N a_n e^{in^2 t}\Bigr\|_{L^2_t([0,2\pi])}\lesssim \| a\|_{\l^2}.
$$
The claim for $2<p<4$ follows by interpolation.  The upper bound we place on $\theta(p)$ is solely there to simplify subsequent analysis.
\end{proof}

\begin{proof}[Proof of Proposition~\ref{P:uN}]
Claim \eqref{uN0} follows from the observation that the Fourier transform of the function $x\mapsto\sech(x)$ is given by
$$
\widehat{\sech}(\xi)=\sqrt{\tfrac\pi2}\sech\bigl(\tfrac{\pi\xi}2\bigr).
$$
In this way, 
$$
\widehat{u_{N,\lambda}}(0,\xi)= (-1)^N\sqrt{\tfrac\pi2}\,N\sech\bigl(\tfrac{\pi\lambda\xi}2\bigr).
$$
Performing the change of variables $\eta=\lambda\xi$ and using that $\lambda\geq 1$, we get
\begin{align*}
\bigl\|u_{N,\lambda}(0)\bigr\|_{\mathcal{F}L^p_s(\R)}^p &\simeq N^p \lambda^{-1} \int_\R \bigl\langle \tfrac\eta\lambda \bigr\rangle^{ps} \sech^p\bigl(\tfrac{\pi\eta}2\bigr)\, d\eta \\
&\lesssim N^p \lambda^{-1} \int_\R \langle \eta\rangle^{p|s|}\sech^p\bigl(\tfrac{\pi\eta}2\bigr)\, d\eta \lesssim N^p \lambda^{-1}.
\end{align*}
Combining this with the lower bound
\begin{align*}
\bigl\|u_{N,\lambda}(0)\bigr\|_{\mathcal{F}L^p_s(\R)}^p \simeq N^p \lambda^{-1} \!\!\int_\R \!\bigl\langle \tfrac\eta\lambda \bigr\rangle^{ps} \sech^p\bigl(\tfrac{\pi\eta}2\bigr)\, d\eta 
\gtrsim N^p \lambda^{-1} \!\!\int_{-1}^1 \!\! \sech^p\bigl(\tfrac{\pi\eta}2\bigr)\, d\eta \gtrsim N^p \lambda^{-1},
\end{align*}
yields \eqref{uN0}.

We now turn to the claims \eqref{uN small} and \eqref{uN big}, for which we will rely on Corollary~\ref{C:asymp}.  Following the notation in Section~\ref{S:3} and accounting for the scaling symmetry of \eqref{mkdv}, we write
$$
v_{N, \lambda} (t,x)= \sum_{j=1}^N \lambda ^{-1} \psi_j\bigl(\lambda^{-3}t,\lambda ^{-1}x\bigr)  
$$
where
$$
\psi_j(t,x)= (-1)^{N}(2j-1)\sech\Bigl((2j-1)\bigl[x-(2j-1)^2t-c_j\bigr]\Bigr).
$$
Correspondingly,
$$
\widehat{\psi_j}(t,\xi)=A_{j}(\xi) e^{-i(2j-1)^2 \xi t} \qtq{where} A_{j}(\xi)=(-1)^{N}\sqrt{\tfrac\pi2}\, e^{-ic_j\xi }\sech\bigl(\tfrac{\pi\xi}{2(2j-1)}\bigr).
$$ 
and so
$$
\widehat{v_{N, \lambda}}(t,\xi)= \sum_{j=1}^NA_{j}(\lambda \xi) e^{-i(2j-1)^2 \lambda^{-2}\xi t}.
$$
Using that the function $\sech$ is radially decreasing, we may simply bound
\begin{align}\label{AJ_Lp}
\bigl\|A_j(\lambda \xi)\bigr\|_{\ell^2_j}^2\simeq \sum_{j=1}^N \sech^2\bigl(\tfrac{\pi\lambda\xi}{2(2j-1)}\bigr)\lesssim N \sech^2\bigl(\tfrac{\pi\lambda \xi}{4N}\bigr).
\end{align}

As an application of Lemma~\ref{L:Bourgain}, we will show that for all $2<p<\infty$ and all $s\in\R$, we have 
\begin{align}\label{FLp psi}
\tfrac{1}{T}\bigl\|\langle \xi\rangle^s\widehat{v_{N, \lambda}}(t,\xi)\bigr\|_{L^p_{t,\xi}([T,2T]\times\R)}^p\lesssim_{p,s} N^{1+\frac{p}{2}+\theta(p)}\lambda^{-1} \qtq{for} \!T\gtrsim \lambda^3N^{\frac p2-1-\theta(p)}.
\end{align}

Performing a change of variables and using Lemma~\ref{L:Bourgain} and \eqref{AJ_Lp}, for $\xi\geq \tfrac{2\pi\lambda^{2}}{T}$ we may bound
\begin{align*}
\tfrac{1}{T}\bigl\|\widehat{v_{N, \lambda}}(\xi)\bigr\|_{L^p_t([T,2T])}^p
&=\tfrac{1}{\lambda^{-2}\xi T}\int_{\lambda^{-2}\xi T}^{2\lambda^{-2}\xi T}\Bigl|\sum_{j=1}^N A_j\bigl(\lambda \xi\bigr)e^{-i(2j-1)^2t}\Bigr|^p\,dt\\
&\le \tfrac{1}{\lambda^{-2}\xi T}\int_{\lambda^{-2}\xi T}^{\lambda^{-2}\xi T+2\pi\bigl\lceil{\frac{\lambda^{-2}\xi T}{2\pi}}\bigr\rceil}\Bigl|\sum_{j=1}^N A_j\bigl(\lambda\xi\bigr)e^{-i(2j-1)^2t}\Bigr|^p\,dt\\
&\lesssim_p \tfrac{1}{\lambda^{-2}\xi T} \bigl\lceil \tfrac{\lambda^{-2}\xi T}{2\pi} \bigr\rceil N^{\theta(p)}\bigl\|A_j\bigl(\lambda\xi\bigr)\bigr\|_{\l^2_j}^p\\
&\lesssim_p   N^{\frac p2+\theta(p)}\sech^p\bigl(\tfrac{\pi\lambda \xi}{4N}\bigr).
\end{align*}
By symmetry, we obtain an analogous bound for $\xi\le -\tfrac{2\pi\lambda^2}{T}$.

On the other hand, for $|\xi|< \frac{2\pi\lambda^2}T$, we simply estimate  
$$
\tfrac{1}{T}\bigl\|\widehat{v_{N,\lambda}}(\xi)\bigr\|_{L^p_{t}([T,2T])}^p\lesssim \Bigl(\sum_{j=1}^N \bigl|A_j\bigl(\lambda\xi\bigr)\bigr|\Bigr)^p\lesssim N^p.
$$

Combining these bounds and using that $\lambda\geq N\geq 1$, we deduce that for $T\gtrsim \lambda^2$,
\begin{align*}
\tfrac{1}{T}\bigl\| &\langle\xi\rangle^s\widehat{v_{N, \lambda}}(t,\xi)\bigr\|_{L^p_{t,\xi}([T,2T]\times\R)}^p\\
&\lesssim \tfrac{1}{T}\int_{|\xi|\ge \frac{2\pi\lambda^2}{T}}\bigl\|\langle\xi\rangle^s\widehat{v_{N, \lambda}}(\xi)\bigr\|_{L^p_{t}([T,2T])}^p\,d\xi+\tfrac{1}{T}\int_{|\xi|\le \frac{2\pi\lambda^2}{T}}\bigl\|\langle\xi\rangle^s\widehat{v_{N, \lambda}}(\xi)\bigr\|_{L^p_{t}([T,2T])}^p\,d\xi\\
&\lesssim_p N^{\frac p2+ \theta(p)}\int_{|\xi|\ge \frac{2\pi\lambda^2}{T}}\langle\xi\rangle^{ps}\sech^p\bigl(\tfrac{\pi\lambda \xi}{4N}\bigr)\, d\xi +N^p \int_{|\xi|\le \frac{2\pi\lambda^2}{T}}\langle\xi\rangle^{ps}\, d\xi\\
&\lesssim_pN^{1+\frac{p}{2}+\theta(p)}\lambda ^{-1} \int_\R \langle\tfrac{N\xi}\lambda\rangle^{p|s|}\sech^p\bigl(\tfrac{\pi\xi}{4}\bigr)\, d\xi+\tfrac{N^p\lambda^2}{T}\\
&\lesssim_{p,s} N^{1+\frac{p}{2}+\theta(p)}\lambda^{-1}+\tfrac{N^p\lambda^2}{T}.
\end{align*}
When $T\gtrsim \lambda^3 N^{\frac p2-1-\theta(p)}$, the first summand dominates, yielding \eqref{FLp psi}.  As $\theta(p)<\frac p2$, this condition also guarantees $T\gtrsim \lambda^2$.

We are now ready to verify \eqref{uN small}, for which $2<p<\infty$.  Using Corollary~\ref{C:asymp} and the scaling symmetry of \eqref{mkdv}, we deduce that 
there exists $t_0(N,\lambda)\in \R$ so that 
\begin{align}\label{5:46}
\bigl\|u_{N,\lambda}(t) - v_{N,\lambda}(t)\bigr\|_{\mathcal{F}L^p_s(\R)}^p \lesssim N^{1+\frac{p}{2}+\theta(p)}\lambda^{-1}  \qtq{for all} t\geq t_0(N,\lambda).
\end{align}
On the other hand, \eqref{FLp psi} guarantees that every interval $[T,2T]$ with $T\gtrsim \lambda^3 N^{\frac p2-1-\theta(p)}$ contains a point $\tau=\tau(T)$ so that
\begin{align}\label{5:47}
\bigl\|v_{N,\lambda}(\tau)\bigr\|_{\mathcal{F}L^p_s(\R)}^p\lesssim_{p,s} N^{1+\frac{p}{2}+\theta(p)}\lambda^{-1}.
\end{align}
Thus, we may find $t_*\gtrsim \max\{t_0(N,\lambda ),  \lambda^3 N^{\frac p2-1-\theta(p)}\}$ so that 
\begin{align}\label{5:48}
\bigl\|u_{N,\lambda}(t_*)\bigr\|_{\mathcal{F}L^p_s(\R)}^p\lesssim_{p,s} N^{1+\frac{p}{2}+\theta(p)}\lambda^{-1},
\end{align}
which proves \eqref{uN small} with $\alpha(p) = \frac12-\frac1p-\frac{\theta(p)}p$.  Note $\alpha(p)>0$ since $\theta(p) < \frac p2 -1$.

We turn now to \eqref{uN big}, where $1\leq p<2$. By \eqref{uN0}, the $L^2$-conservation law for \eqref{mkdv}, and H\"older's inequality,
\begin{align}\label{4:11}
N\lambda^{-\frac12}\simeq
\|u_{N,\lambda}(t)\|_{L^2(\R)}
	&\lesssim\|\jb\xi^s\widehat{u_{N,\lambda}}(t)\|_{L^p}^{\frac{p}{4-p}}\|\jb\xi^{-\frac{ps}{4-2p}}\widehat{u_{N,\lambda}}(t)\|_{L^4}^{\frac{4-2p}{4-p}} 
\end{align}
for any $t\in\R$.  Choosing the $t_*$ from \eqref{uN small} with integrability exponent $4$ and regularity $-\frac{ps}{4-2p}$, we have
$$
\|\jb\xi^{-\frac{ps}{4-2p}}\widehat{u_{N,\lambda}}(t_*)\|_{L^4}\lesssim_{p,s} N^{1-\alpha(4)}\lambda^{-\frac14}
$$
Using this upper bound in \eqref{4:11} and rearranging, we obtain
\begin{align*}
\|u_{N,\lambda}(t_*)\|_{\mathcal F L_s^p}\gtrsim_{p,s} \lambda^{-\frac1p}N^{1 + \alpha(4)\frac{4-2p}{p}}.
\end{align*}
This proves \eqref{uN big} with $\alpha(p) = \alpha(4)\frac{4-2p}{p}$, which is indeed positive.
\end{proof}

\begin{proof}[Proof of Theorem~\ref{T:main}]
We first discuss the case $1\leq p<2$.  For each $N\in\nN$, we choose $\lambda_N= N^{p(1+\frac12\alpha(p))}$, where $\alpha(p)$ is as in Proposition~\ref{P:uN}.

Let $u_{N}(t,x)$ denote the solution to \eqref{mkdv} with initial data
\begin{equation}\label{UNL}
u_{N}(0,x)=(-1)^N \tfrac{N}{\lambda_N}\sech\bigl(\tfrac{x}{\lambda_N}\bigr).
\end{equation}
By \eqref{uN0}, we have
\begin{align*}
\bigl\|u_{N}(0)\bigr\|_{\mathcal{F}L^p_s(\R)}\simeq_{p,s} N \lambda_N^{-\frac1p}\simeq_{p,s} N^{-\frac12\alpha(p)}.
\end{align*}
On the other hand, Proposition~\ref{P:uN} guarantees the existence of a time $t_*$ so that
\begin{align*}
\bigl\|u_{N}(t_*)\bigr\|_{\mathcal{F}L^p_s(\R)} \gtrsim_{p,s} N^{1+\alpha(p)} \lambda_N^{-\frac1p}\gtrsim_{p,s} N^{\frac12\alpha(p)}.
\end{align*}
Thus, the sequence of solutions $\{u_{N}\}_{N\geq 1}$ satisfies the claims of Theorem~\ref{T:main} when $1<p<2$.

Next, we consider $2<p<\infty$.  For each $N\in\nN$, we choose $\lambda=\lambda_N = N^{p(1-\frac12\alpha(p))}$, where $\alpha(p)$ is as in Proposition~\ref{P:uN}. We again write $u_{N}(t,x)$ for the solution to \eqref{mkdv} with initial data \eqref{UNL}.

By Proposition~\ref{P:uN}, there exists a time $t_{*}$ so that
\begin{align*}
\bigl\|u_{N}(t_*)\bigr\|_{\mathcal{F}L^p_s(\R)} \lesssim_{p,s} N^{1-\alpha(p)} \lambda_N^{-\frac1p}\lesssim_{p,s} N^{-\frac12\alpha(p)},
\end{align*}
while 
\begin{align*}
\bigl\|u_{N}(0)\bigr\|_{\mathcal{F}L^p_s(\R)}\simeq_{p,s} N \lambda_N^{-\frac1p}\simeq_{p,s} N^{\frac12\alpha(p)}.
\end{align*}
To conclude, we define
$$
\widetilde{u}_N(t,x) :=u_N\bigl(t_*-t, -x\bigr).
$$
In view of the time-translation and time-reversal symmetries of \eqref{mkdv}, $\widetilde u_N(t,x)$ is also a solution to \eqref{mkdv} and satisfies
\begin{align*}
\bigl\|\widetilde{u}_N(0) \bigr\|_{\mathcal{F}L^p_s(\R)}\lesssim_{p,s} N^{-\frac12\alpha(p)} \qtq{and} \bigl\|\widetilde{u}_N(t_*) \bigr\|_{\mathcal{F}L^p_s(\R)}\simeq_{p,s} N^{\frac12\alpha(p)}.
\end{align*}
Thus, the sequence of solutions $\{\widetilde u_{N}\}_{N\geq 1}$ satisfies the claims of Theorem~\ref{T:main} when $2<p<\infty$.
\end{proof}

\bibliographystyle{siam}
\bibliography{FLbib}

\begin{thebibliography}{10}

\bibitem{Agrawal}
{\sc G.~P. Agrawal}, {\em Nonlinear Fiber Optics}, Academic Press, Oxford, UK,
  {F}ifth~ed., 2013.

\bibitem{MR1029904}
{\sc J.~Bourgain}, {\em On {$\Lambda(p)$}-subsets of squares}, Israel J. Math.,
  67 (1989), pp.~291--311.

\bibitem{MR1820017}
{\sc T.~Cazenave, L.~Vega, and M.~C. Vilela}, {\em A note on the nonlinear
  {S}chr\"odinger equation in weak {$L^p$} spaces}, Commun. Contemp. Math., 3
  (2001), pp.~153--162.

\bibitem{MR2096258}
{\sc A.~Gr\"unrock}, {\em An improved local well-posedness result for the
  modified {K}d{V} equation}, Int. Math. Res. Not.,  (2004), pp.~3287--3308.

\bibitem{MR2529909}
{\sc A.~Gr\"unrock and L.~Vega}, {\em Local well-posedness for the modified
  {K}d{V} equation in almost critical {$\widehat{H^r_s}$}-spaces}, Trans. Amer.
  Math. Soc., 361 (2009), pp.~5681--5694.

\bibitem{HKVZ:1}
{\sc S.~Haque, R.~Killip, M.~Vi\c{s}an, and Y.~Zhang}, {\em Global
  well-posedness and equicontinuity for modified {K}orteweg--de {V}ries
  equations in modulation spaces}, Pure Appl. Anal., 7 (2025), pp.~615--637.

\bibitem{MR4726498}
{\sc B.~Harrop-Griffiths, R.~Killip, and M.~Vi\c{s}an}, {\em Sharp
  well-posedness for the cubic {NLS} and m{K}d{V} in {$H^s(\Bbb R)$}}, Forum
  Math. Pi, 12 (2024), pp.~Paper No. e6, 86.

\bibitem{MR338587}
{\sc R.~Hirota}, {\em Exact envelope-soliton solutions of a nonlinear wave
  equation}, J. Mathematical Phys., 14 (1973), pp.~805--809.

\bibitem{MR3820439}
{\sc R.~Killip, M.~Vi\c{s}an, and X.~Zhang}, {\em Low regularity conservation
  laws for integrable {PDE}}, Geom. Funct. Anal., 28 (2018), pp.~1062--1090.

\bibitem{MR4534495}
{\sc F.~Klaus}, {\em Wellposedness of {NLS} in modulation spaces}, J. Fourier
  Anal. Appl., 29 (2023), pp.~Paper No. 9, 37.

\bibitem{MR3874652}
{\sc H.~Koch and D.~Tataru}, {\em Conserved energies for the cubic nonlinear
  {S}chr\"odinger equation in one dimension}, Duke Math. J., 167 (2018),
  pp.~3207--3313.

\bibitem{MR369963}
{\sc P.~D. Lax}, {\em Periodic solutions of the {K}d{V} equation}, Comm. Pure
  Appl. Math., 28 (1975), pp.~141--188.

\bibitem{MR2319053}
{\sc G.~D. Lyng and P.~D. Miller}, {\em The {$N$}-soliton of the focusing
  nonlinear {S}chr\"odinger equation for {$N$} large}, Comm. Pure Appl. Math.,
  60 (2007), pp.~951--1026.

\bibitem{MR4081534}
{\sc T.~Oh and Y.~Wang}, {\em Global well-posedness of the one-dimensional
  cubic nonlinear {S}chr\"odinger equation in almost critical spaces}, J.
  Differential Equations, 269 (2020), pp.~612--640.

\bibitem{MR4235636}
\leavevmode\vrule height 2pt depth -1.6pt width 23pt, {\em On global
  well-posedness of the modified {K}d{V} equation in modulation spaces},
  Discrete Contin. Dyn. Syst., 41 (2021), pp.~2971--2992.

\bibitem{MR463733}
{\sc J.~Satsuma and N.~Yajima}, {\em Initial value problems of one-dimensional
  self-modulation of nonlinear waves in dispersive media}, Progr. Theoret.
  Phys. Suppl.,  (1974), pp.~284--306.

\bibitem{MR105798}
{\sc S.~Schechter}, {\em On the inversion of certain matrices}, Math. Tables
  Aids Comput., 13 (1959), pp.~73--77.

\end{thebibliography}

\end{document}